\documentclass[12pt]{article}
\usepackage{amssymb,amsmath,latexsym}
\usepackage[colorlinks=true,
linkcolor=webgreen,
filecolor=webbrown,
citecolor=webgreen]{hyperref}

\definecolor{webgreen}{rgb}{0,.5,0}
\definecolor{webbrown}{rgb}{.6,0,0}

\newtheorem{theorem}{Theorem}[section]

\newtheorem{prop}[theorem]{Proposition}

\newtheorem{kor}[theorem]{Corollary}
\newtheorem{coro}[theorem]{Corollary}

\newcommand{\bew}{\paragraph{Proof.}}
\newtheorem{bem}[theorem]{Remark}

\newtheorem{lemma}[theorem]{Lemma}
\newtheorem{proposition}[theorem]{Proposition}
\newtheorem{definition}[theorem]{Definition}

\newcommand{\fa}{\mbox{ for}\mbox{ all }}
\newcommand{\be}{\begin{enumerate}}
\newcommand{\ee}{\end{enumerate}}
\newcommand{\bi}{\begin{itemize}}
\newcommand{\ei}{\end{itemize}}

\newcommand{\ba}{\begin{array}}
\newcommand{\ea}{\end{array}}
\newcommand{\Ra}{\Rightarrow}

\newcommand{\plm}{\mbox{\raisebox{0.4ex}{${\scriptscriptstyle \pm} $}}}
\newcommand{\ra}{\rightarrow}

\newcommand{\zentr}{\mbox{\large{\sf Y}\normalsize}}

\newcommand{\Z}{{\mathbb{Z}}}
\newcommand{\Q}{{\mathbb{Q}}}
\newcommand{\F}{{\mathbb{F}}}
\newcommand{\MM}{{\mathbb{M}}}
\newcommand{\N}{{\mathbb{N}}}
\newcommand{\R}{{\mathbb{R}}}

\newcommand{\RR}{{\mathbb{R}}}
\newcommand{\FF}{{\mathbb{F}}}
\newcommand{\C}{{\mathbb{C}}}

\newcommand{\sC}{{\cal{C}}}
\newcommand{\CC}{{\cal{X}}}
\newcommand{\G}{{\cal{G}}}

\DeclareMathOperator{\wt}{wt}

\DeclareMathOperator{\cwe}{cwe}
\DeclareMathOperator{\fwe}{fwe}
\DeclareMathOperator{\hwe}{hwe}
\DeclareMathOperator{\Lie}{Lie}

\DeclareMathOperator{\End}{End}
\DeclareMathOperator{\diag}{diag}
\DeclareMathOperator{\Aut}{Aut}
\DeclareMathOperator{\Gal}{Gal}

\DeclareMathOperator{\Sym}{Sym}

\newcommand{\0}{{\bf 0}}
\newcommand{\1}{{\bf 1}}

\newcommand{\eb}{\phantom{zzz}\hfill{$\square $}\smallskip}

\renewcommand{\em}{\sf}

\def\trleq{\trianglelefteq}

\begin{document}
\begin{center}
{\Large {\bf  The Invariants of the Clifford Groups}} \\
\vspace{1.5\baselineskip}
{\em Gabriele Nebe}\footnote{Most of this work was carried out during G. Nebe's visit to AT\&T Labs in the Summer of 1999} \\
\vspace*{1\baselineskip}
Abteilung Reine Mathematik der Universit\"at Ulm \\
89069 Ulm, Germany \\
nebe@mathematik.uni-ulm.de \\
\vspace{1\baselineskip}
and \\
\vspace{1\baselineskip}
{\em E. M. Rains} and {\em N. J. A. Sloane} \\
\vspace*{1\baselineskip}
Information Sciences Research, AT\&T Shannon Labs \\
180 Park Avenue, Florham Park, NJ 07932-0971, U.S.A. \\
rains@research.att.com, njas@research.att.com \\
\vspace{1.5\baselineskip}
December 6, 1999; revised September 8, 2000 \\
\vspace{1.5\baselineskip}
{\bf ABSTRACT}
\vspace{.5\baselineskip}
\end{center}

The automorphism group of the Barnes-Wall lattice $L_m$ in dimension $2^m$
$(m \neq 3 )$ is a subgroup of index 2 in a certain ``Clifford group'' $\sC_m$
of structure $2_+^{1+2m} . O^+ (2m,2)$.
This group and its complex analogue $\CC_m$ of structure $(2_+^{1+2m} \zentr Z_8) . Sp(2m,2)$ have
arisen in recent years in connection with the construction of orthogonal
spreads, Kerdock sets, packings in Grassmannian spaces, quantum codes, 
Siegel modular forms and spherical designs.
In this paper we give a simpler proof of Runge's 1996 result
that the space of invariants for $\sC_m$ of degree $2k$ is spanned by the complete weight
enumerators of the codes $C \otimes \FF_{2^m}$, where $C$ ranges over all
binary self-dual codes of length $2k$;
these are a basis if $m \ge k-1$.
We also give new constructions for $L_m$ and $\sC_m$:
let $M$ be the $\Z [ \sqrt{2} ]$-lattice with Gram matrix
$\left[ \begin{array}{rr}
2 & \sqrt{2} \\
\sqrt{2} & 2
\end{array}
\right]$.
Then $L_m$ is the rational part of $M^{\otimes m}$, and $\sC_m = \Aut (M^{\otimes m} )$.
Also, if $C$ is a binary self-dual code not generated by vectors of weight 2,
then
$\sC_m$ is precisely the automorphism group of the complete weight enumerator of $C \otimes \FF_{2^m}$.
There are analogues of all these results for the complex group
$\CC_m$, with ``doubly-even self-dual code'' instead of ``self-dual code''.

\vspace{1\baselineskip}


\clearpage
\setcounter{page}{2}


\section{Introduction}

In 1959 Barnes and Wall \cite{BW59} constructed a family of lattices in dimensions $2^m$, $m=0,1,2, \ldots$.
They distinguished two geometrically similar lattices $L_m \subseteq L'_m$ in $\RR^{2^m}$.
The automorphism group\footnote{More precisely,
$\G_m = \Aut (L_m ) \cap \Aut (L'_m )$ for all $m$, and $\G_m = \Aut (L_m )$ unless $m = 3$.}
$\G_m = \Aut (L_m )$ was investigated in a series of papers by
Bolt, Room and Wall \cite{Bol61}, \cite{BRW61}, \cite{BRW61a}, \cite{Wal62}.
$\G_m$ is a subgroup of index 2 in a certain group ${\cal C}_m$ of structure
$2_+^{1+2m} . O^+ (2m,2)$.
We follow Bolt et~al. in calling ${\cal C}_m$ a Clifford group.
This group and its complex analogue ${\CC}_m$ are the subject of the present paper.

These groups have appeared in several different contexts in recent years.
In 1972 Brou\'{e} and Enguehard \cite{BE73} rediscovered the Barnes-Wall lattices and
also determined their automorphism groups.
In 1995, Calderbank, Cameron, Kantor and Seidel \cite{CCKS97} used
the Clifford groups to construct orthogonal spreads and Kerdock sets,
and asked ``is it possible to say something about [their] Molien series,
such as the minimal degree of an invariant?''.

Around the same time, Runge \cite{Run93}, \cite{Run95}, \cite{Run95a}, \cite{Run95b}
(see also \cite{Duk93}, \cite{Our97})
investigated these groups in connection with Siegel
modular forms.  Among other things, he established
the remarkable result that the space of
homogeneous invariants for ${\cal C}_m$ of degree $2k$ is spanned by the complete
weight enumerators of the codes $C \otimes_{\FF_2} \FF_{2^m}$,
where $C$ ranges over all binary self-dual (or type I) codes of length $2k$,
and the space of homogeneous invariants for $\CC_m$ of degree $8k$ is spanned by the complete
weight enumerators of the codes $C \otimes_{\FF_2} \FF_{2^m}$,
where $C$ ranges over all binary doubly-even self-dual (or type II) codes of length $8k$.
One of our goals is to give a simpler proof of these two assertions,
not involving Siegel modular forms
(see Theorems \ref{theorem218} and \ref{thcxi}).

Around 1996, the Clifford groups also appeared in the study of fault-tolerant
quantum computation and the construction of quantum error-correcting codes
\cite{BDSW96}, \cite{CRSS97}, \cite{CRSS98}, \cite{Kit97}, and in the construction
of optimal packings in Grassmannian spaces \cite{CHRSS99}, \cite{CHS96}, \cite{SS98}.
The story of the astonishing coincidence (involving the 
group ${\cal C}_3$) that led to \cite{CHRSS99}, \cite{CRSS97} and
\cite{CRSS98} is told in \cite{CRSS98}.
(Other recent references that mention these groups are
\cite{Gla95}, \cite{KL88}, \cite{Win72}.)

Independently, and slightly later,
Sidelnikov \cite{Sid97}, \cite{Sid97a}, \cite{Sid99}, \cite{Sid99a} (see also
\cite{Kaz98}) came across the group ${\cal C}_m$ when studying
spherical codes and designs.
In particular, he showed that for $m \geq 3$
the lowest degree harmonic invariant of ${\cal C}_m$
has degree 8, and hence that the orbit under ${\cal C}_m$ of any point on a sphere in
$\RR^{2^m}$ is a spherical 7-design.
(Venkov \cite{Ven00} had earlier shown
that for $m \geq 3$ the minimal vectors of the Barnes-Wall lattices form 7-designs.)

In fact it is an immediate consequence of Runge's results
that for $m \geq 3$ $\sC_m$ has a unique harmonic invariant
of degree 8 and no such
invariant of degree 10 (see Corollary \ref{Mol8}).
The space of homogeneous invariants of degree 8 is spanned by the fourth power of the
quadratic form and the complete weight enumerator of the code
$H_8 \otimes_{\FF_2} \FF_{2^m}$, where $H_8$ is the $[8,4,4]$ Hamming code.
An explicit formula for this complete weight enumerator is given in Theorem \ref{CWE8}.

Our proof of the real version of Runge's theorem
is given in Section 4 (Theorem \ref{theorem218}),
following two preliminary sections dealing with the
group ${\cal C}_m$ and with generalized weight
enumerators.

In Section 5 we study the connection between 
the group ${\cal C}_m$ and the Barnes-Wall lattices.
We define the balanced 
Barnes-Wall lattice $M_m$ to be 
the $\Z[ \sqrt{2} ]$-lattice $\sqrt{2} L'_m + L_m$.
Then $M_m = M_1^{\otimes m}$ (Lemma \ref{tensor}),
which leads to a simple construction: the Barnes-Wall
lattice is just
the rational part of $M_1^{\otimes m}$.
Furthermore $\sC_m = \Aut (M_m )$ (Proposition \ref{tensor2}).
Also, if $C$ is any binary self-dual code that is not generated by vectors of weight 2,
$\sC_m = \Aut (\cwe (C \otimes \FF_{2^m} ))$ (Corollary \ref{tensor3}).
The proof of this makes use of the fact that $\sC_m$ is a maximal finite subgroup of $GL (2^m, \RR )$
(Theorem \ref{maxfin}).
Although there are partial results about the maximality of $\sC_m$ in Kleidman and Liebeck \cite{KL88},
this result appears to be new.
The proof does not use the classification of finite simple groups.

The analogous results for the complex Clifford group
$\CC_m$ are given in Section 6. Theorem \ref{thcxi} is Runge's theorem.
Extending scalars, let
$\MM_m$ be the hermitian $\Z[\zeta _8]$-lattice
$\Z[\zeta _8]\otimes _{\Z [\sqrt{2}]} M_m$.
Then $\CC_m$ is the subgroup of $U(2^m, \Q [ \zeta_8 ])$ preserving
$\MM_m$ (Proposition \ref{tensor4}).
Theorem \ref{maxfin2} shows that, apart from the center, $\CC_m$ is a maximal
finite subgroup of $U(2^m, \C)$, and Corollary \ref{tensor5} is the analogue of
Corollary \ref{tensor3}.

Bolt et al. \cite{Bol61}, \cite{BRW61}, \cite{BRW61a}, \cite{Wal62}
and Sidelnikov \cite{Sid97}, \cite{Sid97a}, \cite{Sid99}
also consider the group $\sC_m^{(p)}$ obtained by replacing 2 in the definition
of $\sC_m$ by an odd prime $p$. In the final section we give
some analogous results for this group.

In recent years many other kinds of self-dual codes have been studied by a number of authors.
Nine such families were named and surveyed in \cite{RS98}.
In a sequel \cite{NRS00a} to the present paper we will give a 
general definition of the ``type'' of a self-dual code which includes all 
these families as well as other self-dual codes over rings and modules.
For each ``type'' we investigate the structure of the associated
``Clifford-Weil group''
(analogous to $\sC_m$ and $\CC_m$ for types I and II) and its ring of invariants.

The results in this paper and in Part II can be regarded as providing a general 
setting for Gleason's theorems \cite{Gle70}, \cite{MS77}, \cite{RS98}
about the weight enumerator of a binary self-dual code (cf. the case
$m=1$ of Theorem \ref{theorem218}), a doubly-even binary self-dual code (cf.
the case $m=1$ of Theorem \ref{thcxi}) and
a self-dual code over $\F _p$ (cf. the case $m=1$ of Theorem \ref{thpi}).
They are also a kind of discrete analogue of a long series of theorems going
back to Eichler (see for example \cite{Boe89},
\cite{Run93}, \cite{Run95}, \cite{Run95b}),
stating that under certain conditions
theta series of quadratic forms are bases for spaces
of modular forms: here complete weight enumerators of
generalized self-dual codes are bases for spaces of invariants of ``Clifford-Weil groups''.

\section{The real Clifford group ${\cal C} _m$}

This initial section defines the real Clifford group ${\cal C} _m$.
The extraspecial $2$-group 
$E(m) \cong 2^{1+2m}_+$ is a subgroup of 
the orthogonal group $O(2^m,\R)$.
If $m=1$ then 
$$E(1):= \left\langle \sigma _1 :=\left(\ba{rr} 0 & 1 \\ 1 & 0 \ea \right) , \sigma _2 :=\left( \ba{rr} 1 & 0 \\
0 & -1 \ea \right) \right\rangle \cong D_8$$
is the automorphism group of the
2-dimensional standard lattice.
In general $E(m)$ is the $m$-fold tensor power of 
$E(1)$:
$$E(m):= E(1)^{\otimes m}  = E(1) \otimes \cdots \otimes E(1)~, $$
and is generated by the tensor products of $\sigma _1 $ and $\sigma _2 $ with 
$2\times 2$ identity matrices $I_2$.

\begin{definition}
The {\em real Clifford group} 
${\cal C} _m $ is  the normalizer in 
$O(2^m,\R)$ of the extraspecial $2$-group  $E(m)$.
\end{definition}

The natural representation of $E(m)$ is absolutely irreducible.
So the centralizer of $E(m)$ in the full orthogonal group
is equal to 
$ \{ \plm I_{2^m} \}$, which is 
the center of $E(m)$.
Then ${\cal C}_m/E(m)$ embeds into 
the outer automorphism group of $E(m)$.
The quotient group 
$E(m)/Z(E(m))$ is isomorphic to a $2m$-dimensional vector space over $\F _2$.
Since every outer automorphism has to respect the $\{+1, -1\}$-valued quadratic form
$$E(m)/Z(E(m)) \cong \F_2^{2m} \ra Z(E(m)) \cong \F_2 ,~~ x\mapsto x^2 ,$$
it follows easily that the outer automorphism group of
$E(m)$ is isomorphic to $O^+(2m,2)$, the full orthogonal group of a quadratic
form of Witt defect 0 over $\F _2$ (see e.g. \cite{Win72}).

Since the group
 $2^{1+2m}_+.O^+(2m,2) $ is a subgroup of $O(2^m,\R )$ (cf. \cite{BRW61a}
or the explicit construction below),
we find that ${\cal C} _m \cong 2^{1+2m}_+.O^+(2m,2)$.
The order of ${\cal C} _m$  is
$$2^{m^2+m+2} (2^m -1) \prod_{j=1}^{m-1} (4^j -1) .$$

To perform explicit calculations we need a convenient set of generators
for ${\cal C}_m$.

\begin{theorem}{\label{generators}}
${\cal C} _m$ is generated by the following elements of $O(2^m,\R ):$
\bi
\item[(1)]
$ \diag ((-1)^{q (v) + a }) $, where $q$ ranges over all $\{0, 1\}$-valued
quadratic forms on $\F _2^{\,m}$ and $a\in \{ 0,1\}$,
\item[(2)] $AGL(m,2)$,
acting on $ \R ^{2^m} = \otimes ^m (\R ^2) =
 \R [\F _{2}^{\,m} ]$ by permuting the
basis vectors in $\F_2^{\,m}$ , and
\item[(3)]
the single matrix
$h \otimes I_2\otimes \cdots \otimes I_2 $
where 
$h:=\frac{1}{\sqrt{2}} \left(\ba{rr} 1 & 1 \\ 1 & -1 \ea \right) $.
\ei
\end{theorem}

\bew
Let $H$ be the group generated by the elements in $(1)$ and  $(2)$.
First, $H$ contains the extraspecial group $E(m)$,
since $ \sigma _1 \otimes I_{2^{m-1}}$ and $ \sigma _2 \otimes I_{2^{m-1}}$ are in $H$
and their images under $GL(m,2)$ generate $E(m)$.

To see that $H/E(m)$ is a maximal parabolic subgroup of $O^+(2m,2)$, note that
by \cite{CCKS97} the elements
 $a\in GL(m,2)$  act on $E(m)/Z(E(m)) \cong \F _2^{2m}$ as
$\left( \ba{cc} a & 0 \\ 0 & a^{-tr} \ea \right)$, and
the elements 
$\diag ((-1)^{q (v) }) $ act as
$\left( \ba{cc} 1 & b \\ 0 & 1 \ea \right)$, where $b$ is the skew-symmetric
matrix corresponding to
the bilinear form $b_q(x,y):=q(x+y)-q(x)-q(y) $.

Since $h\otimes I_{2^{m-1}}  \not\in GL(2^m,\Q )$  is not in $H$,
the group generated by $H$ and this element is ${\cal C}_m$.
\eb

\begin{kor}{\label{cor.2.3}}
${\cal C}_m$ is generated by
$$ \sigma _1 \otimes I_2\otimes \cdots \otimes I_2 , \
 \sigma _2 \otimes I_2\otimes \cdots \otimes I_2 , \
 h\otimes I_2\otimes \cdots \otimes I_2 , \
GL(m,2),\ 
\diag ((-1)^{\Phi (v) }) , $$
where $\Phi $ is the particular quadratic form 
$( \epsilon _1 ,\ldots , \epsilon _m) \mapsto \epsilon _1 \epsilon _2 \in \{ 0,1 \}$
 on $\F_2^{\,m}$.
\end{kor}

\section{Full weight enumerators and complete weight enumerators}

We now introduce certain weight enumerators and show
that they are invariant under the real Clifford group.
Let $C \leq \F _2^N$ be a linear code\footnote{A binary linear code $C$ of
length $N$ is a subspace of $\F _2^N$. If $C \subseteq C^\perp$, $C$ is self-orthogonal;
if $C = C^\perp$, $C$ is self-dual \cite{MS77}, \cite{RS98}.}
of length $N$ over the field $\F _2$.
For $m\in \N$ let $C(m):= C\otimes _{\F _2} \F _{2^m}$ be 
the extension of $C$ to a code over the field with $2^m$ elements.

Let $V$ be the group algebra $V:=\R [\F _{2^m}] = \oplus _{f\in \F _{2^m}} \R x_f$.
Regarding $\F _{2^m} \cong \F_2^{\,m}$ as an $m$-dimensional vector space  over $\F _2$ ,
we have a tensor decomposition 
$$V \cong \otimes ^m (\R ^2)  . $$
In the same manner the group algebra
$\R [C(m)] = \oplus _{c\in C(m)} \R e_c$ embeds naturally into the group algebra
$$\R [ \F_{2^m}^N ] \cong \otimes ^N V = \otimes ^N ( \otimes ^m (\R ^2)) 
\cong \otimes ^m (\otimes ^N (\R ^2)) .$$

\begin{definition}
The {\em full weight enumerator} of $C(m)$ is the element
$$\fwe (C(m)) := \sum _{c \in C(m)} e_c \in \R [C(m)] \subset \otimes ^N V .$$
\end{definition}

(This was called a generalized weight polynomial in \cite{Gle70}
and an exact enumerator in \cite[Chapter 5]{MS77}.)

Fix a basis
 $(a_1,\ldots , a_m)$ of $\F _{2}^{\,m}$ over $\F_2$.
Then a codeword $c\in C(m)$ is just an $m$-tuple 
of codewords in $C$.
The element $c=\sum _{i=1}^m a_ic_i$  corresponds to the 
$m$-tuple $(c_1,\ldots , c_m) \in C^m$, which can also be regarded as
an $m\times N$-matrix $M$ of which the rows are the elements $c_i \in C$.

\begin{lemma}{\label{scalarext}}
Let $$ \fwe _m(C) := \sum _{c_1,\ldots , c_m \in C} 
e_{c_1} \otimes \cdots \otimes e_{c_m} ~ \in ~ \otimes ^m \R [C] ~ \subset ~
\otimes ^m \otimes ^N (\R ^2) .$$
Then the isomorphism $\otimes ^m \R [C] \cong \R [C (m) ]$ 
induced by identifying an $m$-tuple $( c'_1,\ldots , c'_m) \in C^m$
with the codeword $c := \sum _{i=1}^m  a_i c'_i \in C(m)$ maps
$\fwe _m(C)$ onto $\fwe (C(m))$.
\end{lemma}

\bew
Let $c =
 \sum _{i=1}^m  a_i c'_i 
 = (c_1,\ldots , c_N) \in C(m)$. If
$c_i=\sum _{j=1}^m \epsilon ^{(i)}_j a_j $ then
$c'_i = (\epsilon ^{(1)}_i,\ldots , \epsilon ^{(N)}_i) \in C$.
The  generator $e_c$ of $\R [C(m)]$ is 
$$ x_{c_1} \otimes \cdots \otimes x_{c_N} =
(y_{\epsilon ^{(1)}_1} \otimes \cdots \otimes y_{\epsilon ^{(1)}_m})\otimes
\cdots \otimes 
(y_{\epsilon ^{(N)}_1} \otimes \cdots \otimes y_{\epsilon ^{(N)}_m}) \in
\otimes ^N (\otimes ^m (\R ^2)) ,$$
where $\R ^2 = \R [\F _2 ] $ has a basis $y_0,y_1$.
Under the identification above this element is mapped onto
$$(y_{\epsilon ^{(1)}_1} \otimes \cdots \otimes y_{\epsilon ^{(N)}_1})\otimes
\cdots \otimes 
(y_{\epsilon ^{(1)}_m} \otimes \cdots \otimes y_{\epsilon ^{(N)}_m}) \in
\otimes ^m (\otimes ^N (\R ^2)) ,$$
which is the element
$e_{ c'_1}\otimes \cdots \otimes e_{ c'_m} \in
\otimes ^m \R [C]$.
\eb

\begin{definition} (Cf. \cite[Chapter 5]{MS77}.)
The {\em complete weight enumerator} of $C(m)$ is
the following homogeneous polynomial of degree $N$ in $2^m$ variables:
$$\cwe (C(m)) := \sum _{c\in C(m)} \prod _{f \in \F_{2^m}} x_f^{a_f(c)} 
\in \R [x_f \mid f\in \F_{2^m} ] ~,  $$
where $a_f(c)$ is the number of components of $c$ that are equal to $f$.
\end{definition}

\begin{bem}{\label{symm}}
{\rm The complete weight enumerator of $C(m)$ is the projection under 
$\pi $ of the 
full weight enumerator of $C(m)$ to the symmetric power $\Sym_N(V)$,
where $\pi : \otimes ^N V \ra \R [x_f \mid f\in \F _{2^m}]$
is the $\R $-linear mapping defined by
$x_{f_1}\otimes \cdots \otimes x_{f_N}\mapsto 
x_{f_1} \cdots  x_{f_N}:$
 $$\cwe (C(m)) =  \pi (
\fwe (C(m)) ) ~ \in ~ \Sym_N(V).$$ }
\end{bem}

\begin{theorem}\label{cweisinvariant}
Let $C$ be a self-dual code over $\F _2$.
\bi
\item[(i)]
The Clifford group ${\cal C} _m$ preserves the full weight enumerator
$\fwe (C(m))$.
\item[(ii)]
The Clifford group ${\cal C} _m$ preserves the complete weight enumerator
$\cwe (C(m))$.
\ei
\end{theorem}

\bew
Let $N$ be the length of $C$, which is necessarily even.
Then ${\cal C} _m$ acts on $\R [\F _{2^m} ^N] = \otimes ^N (\R ^{2^m}) $
diagonally. 
This action commutes with the projection $\pi : \otimes ^N V \ra \R [x_f
\mid f\in \F_{2^m}]$.
So statement (ii) follows immediately
from (i) by Remark \ref{symm}.
To prove (i) it is enough to consider the generators of ${\cal C} _m$.

The generators 
$ \sigma _1 \otimes  I_2 \otimes \cdots \otimes  I_2  , \
 \sigma _2 \otimes  I_2 \otimes \cdots \otimes  I_2  $ and
 $h\otimes  I_2 \otimes \cdots \otimes  I_2  $
of Corollary \ref{cor.2.3}
 are tensor products of the form 
$x\otimes I_{2^{m-1}}$.
By Lemma \ref{scalarext} it is therefore enough to consider the case $m=1$ 
for these generators.
But then the matrix $ \sigma _1 $ acts as 
$ \sigma _1 \otimes \cdots \otimes  \sigma _1  $ on $\otimes ^N (\R ^2)$, mapping 
a codeword $c = (c_1,\ldots , c_N)\in C$ to $c+\1 =
(c_1+1,\ldots , c_N +1)\in C$,
where $\1 $ is the all-ones vector. Since $C$ is self-dual,
$\1 $ is in $C$ and therefore $ \sigma _1 $ only permutes the
codewords and hence fixes  $\fwe (C)$.
Analogously, the matrix $ \sigma _2 $ changes signs of the components of the codewords
in the full weight enumerator:
if $c=(c_1,\ldots , c_N)$, then $x_{c_i} $ is mapped to 
$(-1)^{c_i} x_{c_i}$.
Since the  codewords in $C$ have even weight,
the tensor product $x_{c_1}\otimes \cdots \otimes x_{c_N}$ is
fixed by $ \sigma _2 \otimes \cdots \otimes  \sigma _2 $.
That $h$ preserves the full weight enumerator follows from
the MacWilliams identity \cite[Chapter 5, Theorem 14]{MS77}.

The generator 
$d:=\diag ((-1)^{\Phi(v) })  = \diag(1,1,1,-1)\otimes I_{2^{m-2}}$
only occurs for $m \geq 2$. By Lemma \ref{scalarext} it
suffices to consider the case $m=2$.
Again by Lemma \ref{scalarext}, we regard $d$ as acting on
pairs $(c,c')$ of codewords in $C$. 
Then $d$ fixes or negates $(x_{c_1}\otimes \cdots \otimes x_{c_N}) \otimes 
(x_{c'_1}\otimes \cdots \otimes x_{c'_N})$, and negates it
if and only if $c$ and $c'$ intersect in an odd number of $1$'s. 
This is impossible since $C$ is self-dual, and so
$d$ also preserves $\fwe (C(m))$.

The remaining generators in $g\in GL(m,2)$ permute the elements of
$\F _{2^m}$.
The codewords $c\in C(m)$ are precisely the elements of the form
$c=\sum _{i=1}^m  a_i  c_i $ with $c_i \in C$ and
$(a_1,\ldots , a_m)$ a fixed $\F _2$-basis for $\F_{2^m}$.
Since $g$ acts linearly on $\F _2^{\,m}$, mapping 
$a_i $ onto $\sum _{j=1}^m g_{ij} a_j$, the word 
$c$ is mapped to 
$\sum _{j=1}^m \sum_{i=1}^m  g_{ij} a_j c_i$ which again is in $C(m)$.
Hence these generators also fix $\fwe (C(m))$.
\eb

\section{The ring of invariants of ${\cal C}_m$}

In this section we establish Runge's theorem
that the complete weight enumerators of the codes  $C(m)$
generate the space of invariants for ${\cal C}_m$.

\begin{definition}
A polynomial $p$ in $2^m$ variables is called a
{\em Clifford invariant of genus $m$} if 
it is  an invariant
for the real Clifford group  ${\cal C}_m$.
Furthermore,
$p$ is called a  {\em parabolic invariant} if it is invariant
under the parabolic subgroup $P$ generated by the elements of type $(1)$ and $(2)$ of Theorem \ref{generators},
 and a {\em diagonal invariant} if it 
is invariant under the group generated by the elements of type $(1)$.
\end{definition}

The following is obvious:

\begin{lemma}
 A polynomial $p$ is a diagonal invariant if and only if
all of its monomials are diagonal invariants.
\end{lemma}

Let $M$ be an $m\times N$ matrix over $\F _2$.  We can associate
a monic monomial ${\mu}_M \in \R [x_f \mid f\in \F_{2^m} ]$
with such a matrix by taking the product of the variables 
associated with its columns.  Clearly all monic monomials are of this form,
and two matrices correspond to the same monic monomial if and only if there
is a column permutation taking one to the other.

\begin{lemma} 
A monic monomial ${\mu}_M$ is a diagonal invariant if and only if
 the rows of $M$ are orthogonal.
\end{lemma}

\bew
It suffices to consider quadratic forms $q_{ij}$ with
$q_{ij}(\epsilon _1,\ldots , \epsilon_m) = \epsilon _i\epsilon _j$ 
($1\leq i \leq j \leq m$);
we easily check that the action of $\diag((-1)^{q_{ij}})$ is to multiply
${\mu}_M$ by $(-1)^k$, where $k$ is the inner product
of rows $i$ and $j$ of $M$; the lemma follows.
\eb

For $g \in GL(m,2) \leq AGL(m,2)$ we have
 $g ({\mu }_M )= {\mu }_{g^{tr}M}$, and  
$b \in \F_2^{\,m} \leq AGL(m,2)$
maps ${\mu }_M$ onto ${\mu }_{M+b}$, where 
the matrix $M+b$ has entries $(M+b)_{ij} = M_{ij} + b_i$.
This implies
that ${\mu}_M$ is equivalent to ${\mu}_{M'}$ under
the action of $AGL(m,2)$ if and only if the binary codes
$\langle M,\1 \rangle$ and $\langle M',\1 \rangle$ are equivalent.
We can thus define a parabolic invariant ${\mu}_m(C)$ 
for any self-orthogonal code
$C$ containing $\1$ and of dimension at most $m+1$ by
\[
\mu_m(C) := \sum_{\substack{M\in \F_2^{{\,m}\times N}\\\langle M,\1\rangle=C}}
\mu_M.
\]
We define ${\mu}_m(C)$ to be $0$ if $\1 \not\in C$ or $\dim (C) >m+1$.  
Since the invariants $\mu_m(C)$ are sums over orbits, we have:

\begin{lemma}  {\label {lemma44}}
 A basis for the space of parabolic invariants of degree $N$ is given by
polynomials of the form ${\mu}_m(C)$ where $C$ ranges over the
equivalence classes of
binary self-orthogonal codes of length $N$ containing $\1 $ and of dimension $\le m+1$.
\end{lemma}

\begin{lemma} {\label {Lem.2.13}}
 For any binary self-orthogonal code $C$ containing $\1 $,
\[
\cwe (C(m)) ~=~ \sum_{\1 \in D\subseteq C}   ~ {\mu}_m(D).
\]
\end{lemma}

\bew From the definition,
\[
\cwe (C(m)) = \sum_M {\mu}_M,
\]
where $M$ ranges over $m\times N$ matrices with all rows in $C$.
Let $M$ be such a matrix.
Then $M$ uniquely determines a subcode 
$D:=\langle M , \1 \rangle $ of $C$; we thus have
\[
\cwe(C(m))
= \sum_{\1\in D\subseteq C} \sum_{\langle M , \1 \rangle = D} \mu_M
= \sum_{\1\in D\subseteq C} \mu_m(D)
\]
as required.
\eb

\begin{theorem}{\label{parabolicbasis}}
A basis for the space of parabolic invariants is given by the
polynomials $\cwe (C(m))$, where $C$ ranges over equivalence classes of
self-orthogonal codes containing $\1 $ and of dimension $\le m+1$.
\end{theorem}

\bew The equations in Lemma \ref{Lem.2.13}
form a triangular system which we can solve for 
the polynomials ${\mu}_m(C)$. In particular,
${\mu}_m(C)$ is a linear combination of the $\cwe (D(m))$ 
for subcodes $\1 \in D \subseteq C$.
\eb

Let $X_P$ denote the
linear transformation
$$x~~\mapsto~~  \frac{1}{|P|} \sum_{g \in P} g \cdot x $$
where $P$ is the parabolic subgroup of  $\sC_m$; that is,
$X_P$ is the operation of averaging over the parabolic subgroup.

\begin{lemma}\label{avglemm}
For any binary self-orthogonal code $C$ of even length $N$ containing $\1 $ and of dimension $N/2-r$,
\begin{align}
X_P (h \otimes I_{2^{m-1}}) \cwe (C(m))
=
\frac{1}{(2^m-1)}
[&
(2^{m-r}-2^r) \cwe (C(m))
\label{eq:avglemm}
\\
&
{}+
2^{-r} \sum_{\substack{C \subset C' \subseteq C^{\prime\perp} \\ [C':C]=2}} \cwe (C'(m))
] .
\notag
\end{align}
\end{lemma}
The final sum is over all self-orthogonal codes $C'$ containing $C$ to
index 2.

\bew
By the MacWilliams identity, we find that
\[
(h \otimes I_{2^{m-1}}) \cwe (C(m))
=
2^{-r}
\sum {\mu}_M,
\]
where $M$ ranges over $m\times N$ matrices such that the first row of $M$
is in $C^\perp$ and the remaining rows are in $C$.  For each code $\1\in
D\subseteq C^\perp$, consider the partial sum over the terms with $\langle
M,\1\rangle=D$.  If $D\subseteq C$, the partial sum is just $\mu_m(D)$,
so in particular is a parabolic invariant.  The other possibility is that
$[D:D\cap C]=2$.  For a matrix $M$ with $\langle M,\1\rangle=D$, define
a vector $v_M\in \F_2^{\,m}$ such that $(v_M)_i=1$ if the $i$th row of $M$
is in $C$, and $(v_M)_i=0$ otherwise.  In particular, the partial
sum we are considering is
\[
\sum_{\substack{\langle M,\1\rangle=D\\v_M=(1,0,0,\dots)}} \mu_M.
\]
If $D$ is not self-orthogonal then this sum is annihilated
by averaging over the diagonal subgroup.  Similarly, if
we apply an element of $AGL_m(2)$ to this sum, this simply has the
effect of changing $v_M$.  Thus, when $D \subseteq D^{\perp}$,
\[
X_P \sum_{\substack{\langle M,\1\rangle=D\\v_M=(1,0,0,\dots)}} \mu_M
=
\frac{1}{|\{v\in \F_2^{\,m}:v\ne 0\}|} \mu_m(D).
\]
Hence
\[
X_P (h \otimes I_{2^{m-1}}) \cwe (C(m))
~=~ 2^{-r} ~
\sum_{\1 \in D\subseteq C} {\mu}_m(D)
+
~ \frac{2^{-r}}{2^m-1}
\sum_{\substack{\1 \in D\subseteq C^\perp \\ [D:D\cap C]=2}}
{\mu}_m(D) ,
\]
where the sums are restricted to self-orthogonal codes $D$.
Introducing a variable $C' = \langle D, C \rangle$ into the second sum
(note that since $D \subseteq C^{\perp}$,
$C' \subseteq C^{\prime\perp}$ precisely
when $D \subseteq D^{\perp}$), this becomes
\[
X_P (h \otimes I_{2^{m-1}}) \cwe (C(m))
=
2^{-r}
\sum_{\1 \in D\subseteq C} {\mu}_m(D)
+
\frac{2^{-r}}{2^m-1}
\sum_{\substack{C\subset C'\subseteq C^{\prime\perp} \\ [C':C]=2}} 
\sum_{\substack{D\subseteq C' \\ D\not \subset C}}
{\mu}_m(D).
\]
Any given $C'$ will, of course, contain each subcode of $C$
exactly once, so we can remove the condition $D\not \subset C$ as follows:
\begin{align*}
X_P (h \otimes I_{2^{m-1}}) &\cwe (C(m))\\
&=
2^{-r}
\sum_{\1 \in D\subset C}
{\mu}_m(D)
+
\frac{2^{-r}}{2^m-1}
\sum_{\substack{C\subset C'\subseteq C^{\prime\perp} \\ [C':C]=2}}
\sum_{\1\in D\subseteq C'}
{\mu}_m(D)
\\
&\phantom{{}={}}{}- (2^{2r}-1)
\frac{2^{-r}}{2^m-1}
\sum_{\1\in D\subseteq C}
{\mu}_m(D)\\
&= 
\frac{1}{2^m-1}
[
(2^{m-r}-2^r) \cwe (C(m))
+
2^{-r} 
\sum_{\substack{C\subset C'\subseteq C^{\prime\perp} \\ [C':C]=2}}
\cwe (C'(m))
],
\end{align*}
as required.
\eb

\begin{lemma} \label{lemma48}
 Let $V$ be a finite dimensional 
vector space, $M$ a linear transformation on $V$,
and $P$ a partially ordered set.  Suppose there exists a spanning set $v_p$ of $V$
indexed by $p \in P$ on which $M$ acts triangularly; that is,
\[
M v_p = \sum_{q\ge p} c_{pq} v_q,
\]
for suitable coefficients $c_{pq}$.  Suppose furthermore that $c_{pp}=1$ if
and only if $p$ is maximal in $P$.  Then the fixed subspace of $M$ in $V$
is spanned by the elements $v_p$ for $p$ maximal.
\end{lemma}

\bew
Since the matrix $C = (c_{pq})$ is triangular, there exists another triangular matrix
$D$ that conjugates $C$ into Jordan canonical form.  Setting
\[
w_p = \sum_{q\ge p} d_{pq} v_q,
\]
($d_{pp}\ne 0$), we find
\[
M w_p = \sum_{q\ge p} c'_{pq} w_q,
\]
with $c'_{pp}=c_{pp}$ and $(M-c_{pp}I)^n w_p = 0$ for sufficiently large
$n$.  In other words, each $w_p$ is in the Jordan block
of $M$ with eigenvalue $c_{pp}$.  But the vectors $w_p$ span $V$;
it follows that the Jordan blocks of $M$ on $V$ are spanned by the
corresponding Jordan blocks of $C$.  In particular, this is true for
the block corresponding to $1$.
\eb

\begin{theorem}\label{theorem218} (Runge \cite{Run95b}.)
Fix integers $k$ and $m\ge 1$.
The space of homogeneous
invariants of degree $2k$ for the Clifford group $\sC_m$ of genus $m$ is spanned by
$\cwe (C(m))$, where $C$ ranges over all binary self-dual codes of length $2k$;
this is a basis if $m\ge k-1$.
\end{theorem}

\bew
Let $p$ be a parabolic invariant.
If $p$ is a Clifford invariant then
\[
X_P (h \otimes I_{2^{m-1}}) p = p.
\]
By Lemma \ref{avglemm}, the operator $X_P(h\otimes I_{2^{m-1}})$ acts
triangularly on the vectors $\cwe_m(C)$ (ordered by inclusion); since
\[
\frac{2^{m-r}-2^r}{2^m-1} = 1 \Ra r=0,
\]
the hypotheses of Lemma \ref{lemma48} are satisfied.
The first claim then follows by Lemma \ref{lemma48} and 
Theorem \ref{cweisinvariant}.
Linear independence for $m \ge k-1$ follows from Lemma \ref{lemma44}.
\eb

In fact a stronger result holds:

\begin{theorem}{\label{thavg}}
For any binary self-orthogonal code $C$ of even length $N$ containing $\1 $
and of dimension $N/2-r$,
\[
\frac{1}{|{\cal C}_m|} \sum_{g\in {\cal C}_m} g \cdot \cwe(C(m))
=
\prod_{1\le i\le r} (2^m+2^i)^{-1}
\sum_{C'} \cwe(C'(m)),
\]
where the sum on the right is over all self-dual codes $C'$ containing $C$.
\end{theorem}

To see that this is indeed stronger than Theorem
\ref{theorem218}, we observe that if $p$ is an invariant
for ${\cal C}_m$ then
\[
\frac{1}{|{\cal C}_m|} \sum_{g\in {\cal C}_m} g \cdot p ~=~ p ~.
\]
Since the space of parabolic invariants contains the
space of invariants, the same is true of the span of
\[
\frac{1}{|{\cal C}_m|} \sum_{g\in {\cal C}_m} g \cdot p
\]
where $p$ ranges over the parabolic invariants. By Theorem
\ref{thavg} each of these can be written as a linear combination
of complete weight enumerators of self-dual codes,
and thus Theorem \ref{theorem218} follows.

\bew
For any self-orthogonal code $C$, let
\[
E_m(C) ~ := ~
\frac{1}{|{\cal C}_m|} \sum_{g\in {\cal C}_m} g \cdot \cwe(C(m)).
\]
Averaging both sides of equation \eqref{eq:avglemm} in Lemma
\ref{avglemm} over ${\cal C}_m$,
we find
\[
E_m(C)
=
\frac{1}{(2^m-1)}
[
(2^{m-r}-2^r) E_m(C)
+
2^{-r} \sum_{\substack{C \subset C' \subseteq C^{\prime\perp} \\ [C':C]=2}}
E_m(C')
],
\]
and solving for $E_m(C)$ gives
\[
E_m(C)
=
\frac{1}{(2^r-1)(2^m+2^r)}
\sum_{\substack{C \subset C' \subseteq C^{\prime\perp} \\ [C':C]=2}} 
E_m(C').
\]
By induction on $r$ (observing that the result follows from Theorem \ref{cweisinvariant} when $r=0$),
we have
\[
E_m(C)
=
\prod_{1\le i\le r} (2^m+2^i)^{-1}
\frac{1}{(2^r-1)}
\sum_{\substack{C \subset C' \subseteq C^{\prime\perp} \\ [C':C]=2}} 
\sum_{C' \subset C''={C''}^\perp}
\cwe_m(C'').
\]
But each code $C''$ is counted $2^r-1$ times (corresponding to the
1-dimensional subspaces of $C''/C$); thus eliminating the sum over $C'$
gives the desired result.
\eb

Note that
\[
\cwe (C(m)) (x_0,\ldots, x_{2^{m-1}-1} ,x_0,\ldots, x_{2^{m-1}-1})
=
\cwe (C(m-1)) (x_0,\ldots, x_{2^{m-1}-1}) ~.
\]
This gives a surjective map from the space of genus
$m$ complete weight enumerators to the space
of genus $m-1$  complete weight enumerators.
By Theorem \ref{theorem218} it follows that this also gives
a surjective map from the genus $m$
invariants to the genus $m-1$ invariants.
(Runge's proof of Theorem \ref{theorem218} 
proceeds by first showing this map is surjective,
using Siegel modular forms, and then arguing that this 
implies Theorem \ref{theorem218}.)
Since by Theorem \ref{parabolicbasis}
the parabolic invariants of degree $N$  become linearly independent
when $m\geq \frac{N}{2} -1 $, we have:

\begin{kor}
Let $\Phi_m(t)$ be the Molien series of the Clifford group of genus $m$.
As $m$ tends to infinity, the series $\Phi_m(t)$ tend monotonically
to
\[
\sum_{k = 0}^{\infty} N_{2k} t^{2k} ~,
\]
where $N_{2k}$ is the number of equivalence classes of self-dual codes of
length $2k$.
\end{kor}
(For the definition of Molien series, see for example
\cite{Ben93} or \cite[Chapter 19]{MS77}.)
\vspace*{+.1in}

Explicit calculations for $m=1,2$ show:

\begin{kor}
 The initial terms of the Molien series of the Clifford group
of genus $m\ge 1$ are given by
\[
1+t^2+t^4+t^6+2 t^8 + 2 t^{10} + O(t^{12}),
\]
where the next term is $2 t^{12}$ for $m=1$, and $3 t^{12}$ for $m>1$.
\end{kor}

Sidelnikov \cite{Sid97a}, \cite{Sid99} showed that the lowest degree of a harmonic invariant
of ${\cal C}_m$ is $8$.
Inspection of the above Molien series gives the following
stronger result.

\begin{kor}{\label{Mol8}}
The smallest degree of a harmonic invariant of ${\cal C}_m$ is $8$, and
there is a unique harmonic invariant of degree $8$.
There are no harmonic invariants of degree $10$.
\end{kor}

The two-dimensional space of
homogeneous invariants for $\sC_m$ of degree 8 is spanned by the fourth power of the quadratic form and by $h_m := \cwe (H_8 \otimes_{\FF_2} \FF_{2^m} )$, where $H_8$ is the $[8,4,4]$ binary Hamming code.
We can give $h_m$ explicitly.
\begin{theorem}\label{CWE8}
Let $G(m,k)$ denote the set of $k$-dimensional subspaces of $\FF_2^{\,m}$.
Then
\begin{eqnarray}\label{Eqh8}
h_m & = & \sum_{v \in \FF_2^{\,m}} x_v^8 +
14 \sum_{U \in G(m,1)}
\sum_{d \in \FF_2^{\,m} /U}
\prod_{v \in d+U} x_v^4 \nonumber \\
&& \quad + 168 \sum_{U \in G(m,2)} \sum_{d \in \FF_2^{\,m} / U}
\prod_{v \in d+U} x_v^2 
+ 1344 \sum_{U \in G(m,3)}
\sum_{d \in \FF_2^{\,m}/U} \prod_{v \in d+U} x_v ~.
\end{eqnarray}
\end{theorem}

The second term on the right-hand side is equal to $14 \sum_{\{u,v\}} x_u^4 x_v^4$, where $\{u,v\}$ runs through unordered pairs of elements of $\FF_2^{\,m}$.
The total number of terms is 
$$2^m + 14 \left[ \begin{array}{c} m\\1 \end{array} \right] 2^{m-1} + 168 \left[ \begin{array}{c} m \\ 2 \end{array}\right] 2^{m-2} + 1344
\left[ \begin{array}{c} m \\ 3 \end{array} \right] 2^{m-3} = 2^{4m} ,$$
where
$\left[ \begin{array}{c} m \\ k \end{array} \right] = |G(m,k)|$.

\bew
We will compute $\cwe (H_8 \otimes \FF_2^{\,m})$ (which is equal to $\cwe (H_8 \otimes \FF_{2^m} )$).
Let $H_8$ be defined by the generator matrix
$$
\begin{array}{cccccccc}
0 & 0 & 0 & 0 & 1 & 1 & 1 & 1 \\
0 & 0 & 1 & 1 & 0 & 0 & 1 & 1 \\
0 & 1 & 0 & 1 & 0 & 1 & 0 & 1 \\
1 & 1 & 1 & 1 & 1 & 1 & 1 & 1
\end{array}
~.
$$
A codeword corresponds to a choice of $(a,b,c,d) \in \FF_2^{\,m}$, one for each row;
from the columns of the generator matrix we find that
the corresponding term of the weight enumerator is
$$x_d x_{c+d} x_{b+d} x_{b+c+d} x_{a+d} x_{a+c+d} x_{a+b+d} x_{a+b+c+d} ~.
$$
This depends only on the affine space $\langle a,b,c \rangle +d$.
The four terms on the right-hand side of Eq. (\ref{Eqh8}) correspond to $\dim \langle a,b,c \rangle =0,1,2,3$;
the coefficients are the number of ways of choosing $(a,b,c,d)$ for a given affine space.
If $\dim \langle a,b,c \rangle =3$, for example, there are
$7 \cdot 6 \cdot 4$ ways to choose $a,b,c$ and 8 ways to choose $d$, giving
the coefficient $8 \cdot 7 \cdot 6 \cdot 4 = 1344$.
\eb

\subsection*{Remarks}
\hspace\parindent
(1) The unique harmonic invariant of degree 8 integrates
to zero over the sphere, and so must have zeros on the sphere.
The orbit of any such point under ${\cal C}_m$ therefore forms
a spherical 11-design, cf. \cite{GS81}. This was
already observed by Sidelnikov \cite{Sid99a}.

(2) The case $m=1$:
$\sC_1$ is a dihedral group of order 16 with Molien series
$1/ (1 - \lambda^2) (1- \lambda^8)$, as in Gleason's
theorem on the weight enumerators of binary self-dual codes
\cite{Gle70}, \cite[Problem 3, p.~602]{MS77},
\cite{RS98}.

(3) The case $m=2$:
$\sC_2$ has order 2304 and Molien series
$$\frac{1+ \lambda^{18}}{(1- \lambda^2) (1- \lambda^8 ) (1- \lambda^{12} ) (1- \lambda^{24} )} ~.$$
(The reflection group $[3,4,3]$, No.~28
on the Shephard-Todd list,
cf. \cite[p. 199]{Ben93}, is a subgroup of $\sC_2$ of index 2.)
The unique harmonic invariants $f_8$ and $f_{12}$ (say) of degrees 8 and 12 are easily computed,
and then one can find real points
$(x_{00}, x_{01}, x_{10}, x_{11} ) \in S^3$ where both $f_8$ and $f_{12}$
vanish.
Any orbit of such a point under $\sC_2$ forms a spherical 15-design of size
2304 (cf. \cite{GS81}).
We conjecture that such points exists for all $m \geq 2$.

(4) The group $\sC_3$ of order 5160960 has appeared in sufficiently many different contexts that it is worth placing its Molien series on record.
It is $p( \lambda ) / q( \lambda )$, where $p( \lambda )$ is the
 symmetric polynomial of degree 154 beginning
\begin{eqnarray*}
&& 1 + \lambda^8 + \lambda^{16} + 2 \lambda^{20} + \lambda^{22} + 2 \lambda^{24} + 3 \lambda^{26} + 4 \lambda^{28} \\
& + &
2 \lambda^{30} + 5 \lambda^{32} + 4 \lambda^{34} + 7 \lambda^{36} + 6 \lambda^{38} + 7 \lambda^{40} \\
& + & 8 \lambda^{42} + 11 \lambda^{44} + 9 \lambda^{46} + 12 \lambda^{48} + 13 \lambda^{50} + 14 \lambda^{52} \\
& + & 15 \lambda^{54} +
17 \lambda^{56} + 17 \lambda^{58} + 20 \lambda^{60} + 19 \lambda^{62} \\
& + & 20 \lambda^{64} + 20 \lambda^{66} + 25 \lambda^{68} + 22 \lambda^{70} + 22 \lambda^{72} \\
& + & 24 \lambda^{74} + 25 \lambda^{76} + \cdots ~,
\end{eqnarray*}
and
$$q( \lambda) = (1- \lambda^2 ) (1- \lambda^{12} ) (1- \lambda^{14} ) (1- \lambda^{16} ) (1- \lambda^{24} )^2
(1- \lambda^{30} ) (1- \lambda^{40} ) ~.
$$

(5) For completeness, we mention that the Molien series for $E(1)$ is
$\frac{1}{(1- \lambda^2 ) (1- \lambda^4 )}$, with basic invariants $x_0^2 + x_1^2$ and $x_0^2 x_1^2$.
For arbitrary $m$ the Molien series for $E(m)$ is
$$\frac{1}{2n^2}
\left\{
\frac{1}{(1- \lambda )^n}
+ \frac{1}{(1+ \lambda)^n} +
\frac{n^2 + n-2}{(1- \lambda^2)^{n/2}} +
\frac{n^2 - n}{(1+ \lambda^2 )^{n/2}}
\right\} ~,
$$
where $n=2^m$.

\section{Real Clifford groups and Barnes-Wall-lattices}

In a series of papers \cite{BW59}, \cite{Bol61}, \cite{BRW61}, \cite{BRW61a},
\cite{Wal62}, Barnes, Bolt, Room and Wall investigated a
family of lattices in $\Q ^{2^m}$ (cf. also \cite{BE73}, \cite{SPLAG}).
They distinguish two geometrically similar lattices
$L_m \subseteq L'_m$ in each dimension $2^m$, for which if $m \neq 3$
the automorphism groups $\Aut(L_m) = \Aut(L'_m)$
are subgroups ${\cal G}_m$
of index 2 in the real Clifford group ${\cal C}_m$.
When $m=3$, $L_3$ and $L'_3$ are two versions of the root lattice $E_8$,
and ${\cal G}_3 := \Aut(L_3) \cap \Aut(L'_3)$ has index $270$
in $\Aut(L_3)$ and index 2 in ${\cal C}_3$.

The lattices $L_m$ and $L'_m$ can be defined in terms of an orthonormal
basis $b_0,\ldots , b_{2^m-1}$ of $\R^{2^m}$ as follows.
Let $V := \F _2^{\,m}$ and index the basis elements 
 $b_0,\ldots , b_{2^m-1}$ by the elements of $V$.
For each affine subspace $U \subseteq V $ let 
$\chi _U\in \Q^{2^m} $ correspond to  the characteristic function of $U$:
$\chi _U := \sum _{i=1}^{2^m} \epsilon _i b_i $,
where $\epsilon _i = 1$ if $i$ corresponds to an element of $U$ and
$\epsilon _i = 0$ otherwise.
Then  $L_m$ (resp. $L'_m$) is spanned by the set
$$ \{ 2^{\lfloor (m-d+\delta )/2 \rfloor} \chi _U \mid 
0\leq d \leq m, U ~ \mbox{is a $d$-dimensional
affine subspace of } V \} ~,$$
where $\delta = 1$ for $L_m$ and $\delta =0 $ for $L'_m$.

Extending scalars, we define the $\Z [\sqrt{2}]$-lattice
$$M_m := \sqrt{2} L'_m +  L_m ~,$$
which we call the {\em balanced Barnes-Wall lattice}.

 From the generating sets for $L_m$ and $L'_m$ we have:

\begin{bem}
{\rm $M_m$ is generated by the vectors
$\sqrt{2}^{m-d} \chi _U$, where $0\leq d \leq m$ and $U$ runs through the
affine subspaces of $V$ of dimension $d$.}
\end{bem}

\begin{lemma}{\label{tensor}}
For all $m > 1$, the lattice $M_m$ is a tensor product:
$$M_m = M_{m-1} \otimes _{\Z [\sqrt{2}]} M_1
= M_1 \otimes_{\Z [ \sqrt{2}]} M_1 \otimes_{\Z[\sqrt{2}]} \cdots
\otimes_{\Z [ \sqrt{2}]} M_1
\quad\mbox{$($with $m$ factors$)$.}
$$
\end{lemma}

\bew
Write $V = \F _2^{\,m}  = V_{m-1} \oplus V_1 $ as the direct sum of
an $(m-1)$-dimensional vector space $V_{m-1}$ and a 1-dimensional space 
$V_1=\langle v \rangle$, and arrange the basis vectors 
so that 
$b_0,\ldots , b_{2^{m-1}-1}$  correspond to the elements in $V_{m-1}$ and
$b_{2^{m-1}},\ldots , b_{2^{m}-1}$   to the elements in $v+V_{m-1}$.

Let $\sqrt{2}^{m-d}\chi _U$ be a generator for $M_m$,
where $U=a+U_0$ for a $d$-dimensional linear subspace $U_0 $ of $V$
and $a = a_{m-1} + a_1 \in V_{m-1} \oplus V_1$.

If $U_0 \leq V_{m-1}$, then 
$$\sqrt{2}^{m-d} \chi _U = (\sqrt{2}^{m-1-d} \chi _{a_{m-1} +U_0 }) \otimes 
\sqrt{2} \chi _{a_1 } \in  M_{m-1} \otimes _{\Z [\sqrt{2}]} M_1.$$
Otherwise $U_{m-1} := U_0 \cap V_{m-1}$ has dimension $d-1$   
and $U_{0} = U_{m-1} \cup (v_{m-1}+v+U_{m-1})$ for some $v_{m-1} \in V_{m-1}$.
If $v_{m-1} \in U_{m-1}$, then 
$$\sqrt{2}^{m-d} \chi _U = (\sqrt{2}^{m-1-(d-1)} \chi _{a_{m-1}+U_{m-1}}) \otimes \chi _{V_1}.$$
If $v_{m-1} \not\in U_{m-1}$ we have the identity
$$\sqrt{2}^{m-d} \chi _U = 
(\sqrt{2}^{m-1-d} \chi _{a_{m-1}+U_{m-1}+\F _2 v_{m-1}}) \otimes \sqrt{2} \chi _{a_1}
+
(\sqrt{2}^{m-1-(d-1)} \chi _{a_{m-1}+v_{m-1} +U_{m-1}}) \otimes \chi _{V_1}
$$
$$
-
\sqrt{2}  (\sqrt{2}^{m-1-(d-1)} \chi _{a_{m-1}+v_{m-1}+U_{m-1}}) \otimes \sqrt{2} \chi _{a_1}~.$$
Hence $M_m \subseteq M_{m-1}\otimes _{\Z [\sqrt{2}]} M_1$.
The other inclusion follows more easily by similar arguments.
\eb

In view of Lemma \ref{tensor}, we have the following simple and apparently new construction for the Barnes-Wall lattice $L_m$.
Namely, $L_m$ is the rational part of the $\Z [ \sqrt{2} ]$-lattice
$M_1^{\otimes m}$, where $M_1$ is the $\Z[ \sqrt{2}]$-lattice with Gram matrix
$\left[ \begin{array}{rr}
2 & \sqrt{2} \\
\sqrt{2} & 2
\end{array}
\right]$.
For more about this construction see \cite{NRS00}.

\begin{proposition}{\label{tensor2}}
For all $m\geq 1$, the automorphism group
$\Aut(M_m)$ (the subgroup of
the orthogonal group $O(2^m, \R )$ that preserves $M_m$) is isomorphic to
$\sC_m$.
\end{proposition}

\bew
Let  $(v_1,\ldots , v_{2^m})$ be a $\Z $-basis for $L'_m$
such that $(2v_1,\ldots , 2v_{2^{m-1}}, v_{2^{m-1}+1},\ldots , v_{2^m})$
is a $\Z $-basis for $L_m$.
Then 
$(\sqrt{2}v_1,\ldots , \sqrt{2}v_{2^{m-1}}, v_{2^{m-1}+1},\ldots , v_{2^m})$
is a $\Z [\sqrt{2} ]$-basis for $M_m = \sqrt{2} L'_m + L_m $.
Hence $M_m$ has a $\Z $-basis
$(\sqrt{2} v_1 ,\ldots , \sqrt{2} v_{2^m},$ $
2v_1,\ldots , 2v_{2^{m-1}},$ $ v_{2^{m-1}+1},\ldots , v_{2^m})$.
Since the scalar products of the $v_i$ are integral, the $\Z $-lattice $M_m$
with respect to $\frac{1}{2}$ the trace form of the 
$\Z [\sqrt{2}]$-valued standard form on $M_m$ is isometric to
$\sqrt{2} L'_m \perp L_m$.
In particular, the automorphism group of the $\Z [\sqrt{2}] $-lattice
$M_m$ is the subgroup of 
$\Aut(
\sqrt{2} L_m \perp L'_m ) \cong {\cal G}_m \wr S_2 $ that commutes with 
the multiplication by $\sqrt{2}$.
Hence 
$\Aut(M_m) $ contains ${\cal G}_m = \Aut(L_m)\cap \Aut(L'_m)$
as a subgroup of index at most two.
Since 
$$h\otimes I_{2^{m-1}} \in \Aut(M_1)\otimes \Aut(M_{m-1}) \subseteq
\Aut(M_m) ,$$
by Lemma \ref{tensor},
$[\Aut(M_m) : {\cal G}_m]= 2$ and so 
$\Aut(M_m)\cong {\cal C}_m$.
\eb

\begin{lemma}{\label{maxorder}}
If $m\geq 2$,  then
the $\Z $-span (denoted $\overline{\Z [{\cal C}_m] }$) of the matrices in 
 ${\cal C}_m$ acting on  the $2^m$-dimensional $\Z [\sqrt{2}]$-lattice
$M_m$ is
  $\Z [\sqrt{2}]^{2^m\times 2^m}$.
\end{lemma}

\bew
We proceed by induction on $m$.
Explicit calculations show  that the lemma is true for $m=2$ and $m=3$.
If $m\geq 4$ then $m-2 \geq 2$ and by induction
$ \overline{\Z [{\cal C}_{m-2}]} =
  \Z [\sqrt{2}]^{2^{m-2}\times 2^{m-2}}$
and $ \overline{\Z [{\cal C}_{2}]} =
  \Z [\sqrt{2}]^{4\times 4}$.
Since $M_m = M_2\otimes _{\Z [\sqrt{2}]} M_{m-2}$ ,
the automorphism group
of  $M_m$ contains ${\cal C}_{2} \otimes {\cal C}_{m-2}$.
Hence 
$$
  \Z [\sqrt{2}]^{2^{m}\times 2^{m}} \supseteq
 \overline{\Z [{\cal C}_{m}]} \supseteq 
 \overline{\Z [{\cal C}_{m-2}]} \otimes_{\Z [\sqrt{2}]}
 \overline{\Z [{\cal C}_{2}]} =
  \Z [\sqrt{2}]^{2^{m}\times 2^{m}} .$$ 
\eb

We now proceed to show that for $m\geq 2$
the  real Clifford group ${\cal C}_m$ is a maximal finite subgroup of 
$GL(2^m,\R )$.
For the investigation of possible normal subgroups of
finite groups containing ${\cal C}_m$, the notion of
a {\em primitive matrix group} plays a central role.
A matrix group $G\leq GL(V)$ is called {\em imprimitive}
if there is a nontrivial decomposition $V = V_1 \oplus \ldots \oplus
V_s$
of $V$ into subspaces which are permuted under the action of $G$.
$G$ is called {\em primitive} if it is not imprimitive.
If $N$ is a normal subgroup of $G$ then $G$ permutes the isotypic
components of $V_{|N}$.
So if $G$ is primitive, the restriction of $V$ to $N$ is isotypic,
i.e. is a multiple of an irreducible representation.
In particular, since the image of an irreducible representation of an
abelian group $N$ is cyclic, all abelian normal subgroups of $G$
are cyclic.

\begin{lemma}{\label{op}}
Let $m\geq 2$.
Let $G$ be a finite group with ${\cal C}_m \leq G \leq GL(2^m,\R )$
and let $p$ be a  prime.
If $p$ is odd, the maximal normal $p$-subgroup of $G$ is trivial.
The maximal normal $2$-subgroup of $G$ is either
$E(m)$ if $G={\cal C}_m$, or $Z(E(m)) = \langle -I_{2^m} \rangle$
if $G>{\cal C}_m$.
\end{lemma}

\bew
We first observe that the only nontrivial normal subgroup of ${\cal C}_m$
that is properly contained in $E(m)$ is $Z(E(m)) = \langle -I_{2^m} \rangle $.
Therefore, if $U$ is a normal subgroup of $G$, $U\cap E(m)$ is
one of $1$, $Z(E(m))$ or $E(m)$.

The matrix group ${\cal C}_m$ and hence also  $G$ is primitive.
In particular, all abelian normal subgroups of $G$ are cyclic.
Let $p>2$ be a rational prime
and $U \trleq G$ a normal $p$-subgroup of $G$.
The degree of the absolutely irreducible representations
of $U$  that occur in $\R ^{2^m} _{|U}$ is a power of $p$
and divides $2^m$.  So this degree is 1 and $U$ is abelian,
hence cyclic by the primitivity of $G$.
Therefore the automorphism group of $U$ does not contain $E(m)/Z(E(m))$.
Since $C_G(U)\cap E(m)$ is a normal subgroup of ${\cal C}_m$,
it equals $E(m)$ and hence $E(m)$ centralizes $U$.
Since $E(m)$ is already absolutely irreducible, $U$ consists
of scalar matrices in $GL(2^m, \R )$, and therefore $U=1$.
If $p=2$ and $G\neq {\cal C}_m$,
then $U\neq E(m)$, because ${\cal C}_m $ is the largest finite
subgroup of $GL(2^m,\R )$ that normalizes $E(m)$.
Since the normal 2-subgroups of $G$ do not contain an
abelian noncyclic characteristic subgroup, 
the possible normal 2-subgroups are classified
in a theorem of P. Hall (cf. \cite[p. 357]{Hup67}).
In particular they do not contain ${\cal C}_m/Z(E(m))$ 
as a subgroup of their automorphism groups, so again 
$U$ commutes with $E(m)$, and therefore consists only of scalar matrices.
\eb

\begin{theorem}{\label{maxfin}}
Let $m\geq 2$.
Then the real Clifford group ${\cal C}_m$ is a maximal finite subgroup of 
$GL(2^m,\R )$.
\end{theorem}

\bew
Let $G$ be a finite subgroup of  $ GL(2^m,\R )$ 
that properly contains ${\cal C}_m$.
By Lemma \ref{op}, all normal $p$-subgroups of $G$ are central.
By a theorem of Brauer, every representation of a finite group is
realizable over a cyclotomic number field (cf. \cite[\S 12.3]{Ser77}).
In fact, since the natural representation of $G$ is real, it is even true that
$G$ is conjugate to a subgroup of $GL(2^m,K)$ for some 
totally real abelian number field $K$ containing $\Q [\sqrt{2}]$
(cf. \cite[Proposition 5.6]{Dre75}).
Let $K$ be a minimal such field and assume that $G\leq GL(2^m,K)$.
Let $R$ be the ring of integers of $K$.
Then $G$ fixes an $R{\cal C}_m$-lattice.
By Lemma \ref{maxorder} all $R{\cal C}_m$-lattices are of the
 form $I \otimes _{\Z [\sqrt{2}]} M_m$
for some   fractional ideal $I$  of $R$,
the group $G$ fixes all $R{\cal C}_m$-lattices and hence also 
$R\otimes _{\Z [\sqrt{2}]} M_m$.
So any choice of an $R$-basis for $M_m$ gives rise to an
embedding $G \hookrightarrow GL(2^m ,R)$,
by which we may regard $G$ as a group of matrices.
Without loss of generality we may assume that $G=\Aut(R\otimes _{\Z [\sqrt{2}]}
M_m )$.
Then the Galois group $\Gamma := \Gal(K/\Q [\sqrt{2}])$ acts on 
$G$ by acting componentwise on the matrices.
Seeking a contradiction, we assume $K \neq \Q [\sqrt{2}]$.
It is enough 
to show that there is a nontrivial element  $\sigma \in \Gamma $
that acts trivially on $G$, because then 
the matrices in $G$ have their entries in the fixed field of
$\sigma$, contradicting the minimality of $K$.

Assume first that there is an odd prime $p$ ramified in $K/\Q $,
and let $\wp $ be a prime ideal of $R$ that lies over $p$.
Then $p$ is also ramified in $K /\Q [\sqrt{2}]$ and therefore
the action of the ramification group, the stabilizer in 
$\Gamma $ of $\wp $, on $R/\wp $ is not faithful, hence
the first inertia group
 $$\Gamma _{\wp } := \{ \sigma \in \Gal(K/\Q [\sqrt{2}]) 
\mid  \sigma(x) \equiv x \pmod{\wp } \fa x\in R \} $$
is nontrivial (see e.g. \cite[Corollary III.4.2]{FT91}).
Since $G_{\wp } := \{ g\in G \mid g \equiv I_{2^m} \pmod{\wp } \}$
is a normal $p$-subgroup of $G$, $G_{\wp } = 1$
by Lemma \ref{op}.
Therefore all the elements in $\Gamma _{\wp } $ act trivially on $G$,
which is what we were seeking to prove.

So $2$ is the only ramified prime in $K$, which implies that
$K = \Q [\zeta _{2^a} + \zeta _{2^a}^{-1}]$ for some $a\geq 3$,
where $\zeta_t = \exp (2 \pi i /t )$.
If $a=3$, then $K=\Q [\sqrt{2}]$, $G = \Aut(M_m) = {\cal C}_m$ 
and we are done.
So assume $a>3$ and let $\wp $ be the prime ideal
of $R$ over $2$ (generated by $(1-\zeta _{2^a})(1-\zeta _{2^a}^{-1})$)
and let $\sigma \in \Gamma $ be the Galois automorphism
defined by $\sigma(\zeta _{2^a} + \zeta _{2^a}^{-1} ) =
\zeta _{2^a}^{2^{a-1}+1} + \zeta _{2^a}^{-2^{a-1}-1}
= - (\zeta _{2^a} + \zeta _{2^a}^{-1} )$.
Then $id=\sigma ^2 \neq \sigma $ and
$$
(\zeta _{2^a} + \zeta _{2^a}^{-1} )-
\sigma(\zeta _{2^a} + \zeta _{2^a}^{-1} )  =
2(\zeta _{2^a} + \zeta _{2^a}^{-1} ) \in 2\wp .$$
Therefore $\sigma \in \Gamma _{2\wp }$.
Since the subgroup 
$G_{2\wp } := \{ g\in G \mid g \equiv I_{2^m} \pmod{2\wp } \}$
of $G$ is trivial (cf.  \cite[Hilfssatz 1]{Bar78})
one concludes that $\sigma $ acts trivially on $G$,
and thus $G$ is in fact defined over 
$\Q [\zeta _{2^{a-1}} + \zeta _{2^{a-1}}^{-1}]$.
The theorem follows by induction.
\eb

\begin{kor}{\label{tensor3}}
Let $m \geq 1$
and let $C$ be a self-dual code over $\F _2$ that is not generated
by vectors of weight $2$.
Then $${\cal C}_m = \Aut_{O(2^m,\R )}(\cwe (C(m)) .$$
\end{kor}

\bew
The proof for the case $m = 1$ will be postponed to Section 6.
Assume $m \geq 2$.
We first show that the parabolic subgroup
$H\leq {\cal C}_m$ acts irreducibly on the Lie algebra $\Lie(O(2^m,\R ))$,
the set of real $2^m \times 2^m$ matrices $X$ such that
$X = - X^{tr}$.
The group $AGL(m,2)$ acts 2-transitively on
our standard basis $b_0,\ldots , b_{2^m-1}$ for $\R^{2^m}$.
A basis for $\Lie(O(2^m,\R ))$ is given by the
matrices
$b_{ij} := b_i\otimes b_j - b_j \otimes b_i$ for $0\leq i < j \leq 2^m-1$.
Since $AGL(m,2)$ acts transitively on the $b_{ij}$,
a basis for the  endomorphism ring $\End_{AGL(m,2)}(\Lie(O(2^m,\R )))$
is given by the orbits of the stabilizer of $b_{01}$.
Representatives for these orbits are 
$b_{01}$, $b_{02}$, $b_{23}$ and $b_{24}$.
But the generator corresponding to the quadratic form
$q(v_1,\ldots , v_m) := v_2^2$ 
negates $b_2$ and fixes $b_0$ and $b_4$, and therefore does not
commute with the endomorphism corresponding to $b_{02}$ or $b_{24}$.
Similarly the endomorphism corresponding to $b_{23}$ 
is ruled out by 
$q(v_1,\ldots , v_m) := v_1v_2$.

Let $G:=\Aut_{O(2^m,\R )}(\cwe(C(m)) $.
Then $G$ is a closed subgroup of $O(2^m,\R )$ 
and hence is a Lie group (cf. \cite[Theorem 3.4]{PR94}).
Since $G$ contains ${\cal C}_m$
it acts irreducibly on $\Lie(O(2^m,\R ))$.
Assume that $G \neq {\cal C}_m$.
Then $G$ is infinite by  Theorem \ref{maxfin} and 
therefore $G$ contains $SO(2^m,\R )$.
However, the ring of invariants of $SO(2^m,\R )$
is generated by the quadratic form
$\sum_{i=0}^{2^m -1} x_{b_i}^2$.
The only binary self-dual codes $C$ that produce such complete weight enumerators
are direct sums of copies of the code $\{00, 11\}$.
\eb

\section{The complex Clifford groups and doubly-even codes}
There are analogues for the complex Clifford group $\CC_m$ for most of the above results.
($Z_a$ will denote a cyclic group of order $a$.)

\begin{definition}
The {\em complex Clifford group} $\CC _m$ is the normalizer in
$U(2^m,\Q [\zeta _8])$ of
the central product $E(m) \zentr Z_4$.
\end{definition}

As in the real case, one concludes that 
\[
\CC _m \cong (2_+^{1+2m} \zentr Z_8 ). Sp(2m,2)
\cong
(2_+^{1+2m} \zentr Z_8 ). O(2m+1,2)
\]
(cf. \cite[Cor. 8.4]{Neb98}).

The analogue of Theorem \ref{theorem218} is the following,
which can be proved in essentially the same way.

\begin{theorem} {\label {thcxi}} (Runge \cite{Run95b}.)
Fix integers $N$ and $m\ge 1$.  The space of homogeneous
invariants of degree $N$ for the complex Clifford group ${\CC}_m$ is spanned by
$\cwe (C(m))$, where $C$ ranges over all binary doubly-even
self-dual codes of length $N$. (In particular, when $N$ 
is not a multiple of $8$, the invariant space is empty.)
\end{theorem}

The analogues of Theorem \ref{thavg} and  Proposition \ref{tensor2} are:

\begin{theorem}
For any doubly-even binary code $C$ of length $N\equiv 0(8)$ containing
$\1$ and of dimension $N/2-r$,
\[
\frac{1}{|{\CC}_m|} \sum_{g\in {\CC}_m} g \cdot \cwe(C(m))
=
\prod_{0\le i<r} (2^m+2^i)^{-1}
\sum_{C'} \cwe(C'(m)),
\]
where the sum is over all doubly-even self-dual codes $C'$ containing $C$.
\end{theorem}

\begin{prop}\label{tensor4}
Let $\MM_m := \Z[\zeta_8] \otimes _{\Z [\sqrt{2}]} M_m$.
Then the subgroup of
$U(2^m, \Q [ \zeta_8 ])$ preserving $\MM_m$ is
precisely $\CC_m$.
\end{prop}

We omit the proofs.

For the analogue of Lemma \ref{maxorder}, observe that
the matrices in  $\CC _m$ generate a maximal order.
Even for $m=1$ the $\Z $-span of the matrices in $\CC _1$ acting on $\MM _1$ 
is the maximal order $\Z [\zeta _8] ^{2\times 2}$.
Hence the induction argument used to prove Lemma \ref{maxorder}
shows that 
$\overline{\Z [\CC _m ]} = \Z [\zeta _8]^{2^m\times 2^m}$.
Therefore the analogue of Theorem \ref{maxfin}  holds even for $m=1$:

\begin{theorem}\label{maxfin2}
Let $m \geq 1$ and 
let $G$ be a finite group such that
$\CC _m \leq G \leq U(2^m,\C )$.
Then 
there exists a root of unity $\zeta$ such that
\[ 
G = \langle \CC _m,\zeta I_{2^m}\rangle.
\]
\end{theorem}

\bew
As in the proof of Theorem \ref{maxfin}, we may assume that
$G$ is contained in $U(2^m,K)$ for some abelian number field $K$ 
containing $\zeta _8$.
Let $R$ be the ring of integers in $K$ and
$T$ the group of roots of unity in $R$.
Then $T \CC _m $ is the normalizer in $U(2^m,K)$ of
$T E(m)$ (cf. \cite[Cor. 8.4]{Neb98}).
As before, the 
 $R\CC_m $-lattices in the
natural module are of the form $I\otimes _{\Z [\zeta _8]} \MM_m$,
where $I$ is a fractional ideal of $R$.
Since $G$ fixes one of these lattices, it also fixes
 $R\otimes _{\Z [\zeta _8]} \MM_m$.
As in the proof of Theorem \ref{maxfin}, we write the elements of $G$ as matrices with
respect to a basis for $\MM_m$ and assume that $G$ is the full 
(unitary) automorphism group
of $R\otimes_{\Z [\zeta _8]} \MM_m$.
Then the  Galois group $\Gamma := \Gal(K/\Q [\zeta _8])$ acts on 
$G$.
 Assume that $G \neq T \CC _m$.
Then $T E(m)$ is not normal in $G$.
As in Lemma \ref{op} one shows that
the maximal normal $p$-subgroup of $G$
is central for all primes $p$.
Let $\wp $ be a prime ideal in $R$ that ramifies in $K/\Q[\zeta _8]$,
and let $\sigma $ be an element of 
 the inertia group $\Gamma _{\wp }$.
Then for all $g\in G$, the image $g^{\sigma }$ satisfies
$a(g):= g^{-1} g^{\sigma } \in 
G_{\wp } := \{ g\in G \mid g \equiv I_{2^m} \pmod{\wp } \}$.
Since $G_{\wp }$ is a normal $p$-subgroup,
where $p$ is the rational prime divisible by $\wp$,
it is central.
Therefore the map $ g \mapsto a(g)$ is a homomorphism of $G$ into
an abelian group, and hence the commutator subgroup $G'$ is fixed under
$\sigma $. 
Since any abelian extension $K$ of $\Q $ that properly contains $\Q [\zeta _8]$
is ramified at some finite prime of $\Q [\zeta _8]$,
we conclude that 
$G ' \subseteq \Aut(\MM_m)$.
Since $E(m)\zentr Z_8 \leq \Aut (\MM _m)' \zentr Z_8$ 
is characteristic in $\Aut(\MM_m)$ and therefore also in $G'\zentr Z_8$,
the group $TE(m)$ is normal in $G$, which is a contradiction.
\eb

\begin{coro}
\label{tensor5}
Assume $m \geq 1$ and
let $C$ be a binary self-dual doubly-even code of length $N$.
Then
$$\Aut_{U(2^m, \C )}
(\cwe (C \otimes \FF_{2^m} )) = \langle \CC_m, \zeta_N I_{2^m} \rangle ~.
$$
\end{coro}

\subsection*{Remarks}
\hspace\parindent
(1) The case $m=1$:
$\CC_1$ is a unitary reflection group (No. 9 on the Shephard-Todd list)
of order 192 with Molien series
$1 /(1- \lambda^8) (1- \lambda^{24} )$, as in Gleason's theorem on the weight
enumerators of doubly-even binary self-dual codes
\cite{Gle70}, \cite[p. 602, Theorem 3c]{MS77}, \cite{RS98}.

(2) The case $m=2$:
$\CC_2$ has order 92160 and Molien series
$$\frac{1+ \lambda^{32}}{(1- \lambda^8) (1- \lambda^{24})^2 (1- \lambda^{40} )} ~.
$$
This has a reflection subgroup of index 2,
No.~31 on the Shephard-Todd list.

(3) The case $m=3$:
$\CC_3$ has order 743178240, and the Molien series can be written as
$p( \lambda^8 ) /q ( \lambda )$, where $p( \lambda)$ is the symmetric
polynomial of degree 44 beginning
\begin{align*}
1+{}&\lambda^{3}+3\lambda^{4}+3\lambda^{5}+6\lambda^{6}+8\lambda^{7}+12\lambda^{8}+18\lambda^{9}+25\lambda^{10}+29\lambda^{11}+40\lambda^{12}+50\lambda^{13}\\
{}+{}&58\lambda^{14}+69\lambda^{15}+80\lambda^{16}+85\lambda^{17}+96\lambda^{18}+104\lambda^{19}+107\lambda^{20}+109\lambda^{21}+112\lambda^{22}+\dots
\end{align*}
and
\[
q( \lambda ) = (1-\lambda^8)(1-\lambda^{16})(1-\lambda^{24})^2(1-\lambda^{40})(1-\lambda^{56})(1-\lambda^{72})(1-\lambda^{120}).
\]
Runge \cite{Run95} gives the Molien
series for the commutator subgroup ${\cal H}_3 = \CC'_3$, of index 2 in $\CC_3$.
The Molien series for $\CC_3$ consists of the terms in the series
for ${\cal H}_3$ that have exponents divisible by 4.
Oura \cite{Our97} has computed the Molien series for ${\cal H}_4 = \CC'_4$,
and that for $\CC_4$ can be obtained from it in the same way.
Other related Molien series can be found in \cite{BDHO99}.

{\bf Proof of Corollary \ref{tensor3}, case $m = 1$.}

Let $C$ be a self-dual binary code of length $n$ with Hamming weight
enumerator $\hwe_C(x,y)$. We will show that if $C$ is not generated
by vectors of weight 2 then
$\Aut_{O(2)} (\hwe_C) = {\cal C}_1$.

Certainly $G:=\Aut_{O(2)}(\hwe_C)$ contains ${\cal C}_1=D_{16}$; we must
show it is no larger.  The only closed subgroups of $O(2)$
containing $D_{16}$ are the dihedral groups $D_{16k}$ for $k\ge 1$
and $O(2)$ itself. 
So if the result is false then $G$ contains a rotation 
\[
\rho(\theta)=
\pmatrix
\cos(\theta)&\sin(\theta)\\
-\sin(\theta)&\cos(\theta)
\endpmatrix
\]
where $\theta$ is not a multiple of $\pi/4$.

Consider the shadow $S(C)$ of $C$ \cite{RS98}; that is, the
set of vectors $v\in \F_2^n$ such that
\[
\wt(v+w)\equiv \wt(v)\pmod{4},\ \fa w\in C.
\]
The weight enumerator of $S(C)$ is given by
$S(x,y) = 2^{-n/2} \hwe_C(x+y,i(x-y))$.
Then $\rho(\theta)\in G$ if and only if
$S(x,y) = S(e^{i\theta} x,e^{-i\theta} y)$,
or in other words if and only if for all $v\in S(C)$,
$(n-2\wt(v))\theta$ is a multiple of $2\pi$.

Now, pick a vector $v_0\in S(C)$, and consider the polynomial $W(x,y,z,w)$ given by
\[
\sum_{v\in C}
~ x^{n-\wt(v_0)-\wt((\1+v_0)\cap v)}
~ y^{\wt((\1+v_0)\cap v)}
~ z^{\wt(v_0)-\wt(v_0\cap v)}
~ w^{\wt(v_0\cap v)} ~.
\]
This has the following symmetries:
$$
W(x,iy,z,-iw) = W(x,y,z,w),
$$
$$
W((x+y)/\sqrt{2},(x-y)/\sqrt{2},(z+w)/\sqrt{2},(z-w)/\sqrt{2}) = W(x,y,z,w).
$$
Furthermore, since $S(C)=v_0+C$, $\rho(\theta)\in G$ if and only if
\[
W(e^{i\theta} x,e^{-i\theta} y,e^{-i\theta} z,e^{i\theta} w)=W(x,y,z,w).
\]
To each of these symmetries we associate a $2\times 2$ unitary
matrix $U$ such that $(x,y)$ is transformed according to $U$ and
$(z,w)$ according to $\overline{U}$.  The first two
symmetries generate the complex group $\CC_1$, which is
maximally finite in $PU(2)$ by Theorem \ref{maxfin2}.  On the other hand,
we can check directly that
\[
\pmatrix e^{i\theta}&0\\0&e^{-i\theta}\endpmatrix
\notin \CC_1,
\]
even up to scalar multiplication.
Thus the three symmetries topologically
generate $PU(2)$; and hence $W$ is invariant under any unitary matrix
of determinant $\pm 1$.  Since
$\hwe_C(x,y) = W(x,y,x,y)$,
it follows that $G=O(2)$.  But then
\[
\hwe_C(x,y) = (x^2+y^2)^{n/2},
\]
implying that $C$ is generated by vectors of weight 2.

This completes the proof of Corollary \ref{tensor3}.
\eb

\section{Clifford groups for $p>2$}

Given an odd prime $p$, there again is a natural representation of
the extraspecial $p$-group $E_p(m) \cong p^{1+2m}_+$ of exponent $p$, this time in $U(p^m,\C)$;
to be precise, $E_p(1)$ is generated by transforms
\[
X:v_x\mapsto v_{x+1},\ \text{and}\ Z:v_x\mapsto \exp(2\pi i x/p) v_x, \
x\in \Z/p\Z  ~,
\]
and $E_p(m)$ is the $m$-th tensor power of $E_p(1)$.  The Clifford
group ${\cal C}^{(p)}_m$ is then defined to be the normalizer in
$U(p^m,\Q[\zeta_{ap}])$ of $E_p(m)$, where $a = gcd \{p+1,4\}$.
As above, one finds that
\[
{\cal C}^{(p)}_m \cong Z_a \times p_+^{1+2m} . Sp(2m,p)
\]
(cf. e.g. \cite{Win72}).

As before, the invariants of these Clifford groups are given by codes:

\begin{theorem}\label{thpi}
Fix integers $N$ and $m\ge 1$.  The space of invariants of degree $N$ for
the Clifford group ${\cal C}^{(p)}_m$ is spanned by $\cwe (C(m))$, where
$C$ ranges over all self-dual codes over $\F_p$ of length $N$ containing
$\1$.
\end{theorem}

\begin{theorem}
For any self-orthogonal code $C$ over $\F_p$ of length $N$ containing
$\1$ and of dimension $N/2-r$,
\[
\frac{1}{|{\cal C}^{(p)}_m|} \sum_{g\in {\cal C}^{(p)}_m} g \cdot \cwe(C(m))
=
\prod_{0\le i<r} (p^m+p^i)^{-1}
\sum_{C'} \cwe(C'(m)),
\]
where the sum is over all self-dual codes $C'$ containing $C$
(and in particular is $0$ if no such code exists).
\end{theorem}

Regarding maximal finiteness, the arguments we used for $p=2$ 
to prove Theorem \ref{maxfin} do not
carry over to odd primes, since the groups ${\cal C}^{(p)}_m$
do not span a maximal order.
Lindsey \cite{Lin70} showed by group theoretic arguments that 
${\cal C}^{(p)}_1$ is a maximal finite subgroup of
$SL(p,\C )$ (cf. \cite{Bli17} for $p=3$,
\cite{Bra67} for $p=5$).
For $p^m = 9$, the theorem below follows from 
\cite{Fei75} and \cite{HW78}.  

\begin{theorem}
Let $p>2$ be a prime and $m\ge 1$.
If $G$ is
a finite group with ${\cal C}^{(p)}_m\leq G\leq GL(p^m,\C)$,
there exists a root of unity $\zeta$ such that
\[ 
G = \langle {\cal C}^{(p)}_m,\zeta I_{p^m}\rangle.
\]
\end{theorem}

\bew
As before we may assume that $G$ is contained in $U(p^m,K)$ for some
abelian number field $K$ containing $\zeta_p$.
Let ${\cal L}$ denote the
set of rational primes $l$ satisfying the following four properties:
(i) $G$ is $l$-adically integral,
(ii) $l$ is unramified in $K$,
(iii) $|G|<|PGL(p^m,l)|$,
(iv) $l$ splits completely in $K$.
Since all but finitely many primes satisfy conditions (i)-(iii), and
infinitely many primes satisfy (iv) (by the {\v C}ebotarev Density Theorem),
it follows that the set ${\cal L}$ is infinite.

Fix a prime ${\mathfrak l}$ over $l\in {\cal L}$.  Since $G$ is
$l$-adically integral, we can reduce it mod ${\mathfrak l}$, obtaining a
representation of $G$ in $GL(p^m,l)$.  Since $p$ is ramified in $K$, $l\ne
p$, so this representation is faithful on the extraspecial group.  Since
the extraspecial group acts irreducibly, the representation is in fact
faithful on the entire Clifford group.  Thus $G\bmod {\mathfrak l}$ contains
the normalizer of an extraspecial group, but modulo scalars is strictly
contained in $PGL(p^m,l)$ (by condition (iii)).
It follows from the main theorem of \cite{KL88} that for $p^m \ge 13$
$G\bmod {\mathfrak l}$ and ${\cal C}^{(p)}_m \bmod {\mathfrak l}$
coincide as subgroups of $PGL(p^m,l)$.
For $p^m < 13$ this already follows from the references in the paragraph
preceding the theorem.

Fix a coset $S$ of ${\cal C}^{(p)}_m$ in $G$.  For each prime ${\mathfrak
l}|l$ with $l\in {\cal L}$, the above argument implies that we can choose
an element $g\in S$ such that $g\propto 1\pmod{\mathfrak l}$.  As there
are infinitely many such primes, at least one such $g$ must get chosen
infinitely often.  But then we must actually have $g\propto 1$ in $K$,
and since $g$ has finite order, $g=\zeta_S$ for some root of unity
$\zeta_S$.

Since this holds for all cosets $S$, $G$ is generated 
by ${\cal C}^{(p)}_m$ together with the roots of unity $\zeta_S$, proving
the theorem.
\eb

\begin{bem}
{\rm
It is worth pointing out that the proof of the main theorem
in \cite{KL88}
relies heavily on the classification of finite simple groups, which
is why we preferred to use our alternative arguments when proving Theorem \ref{maxfin}.
}
\end{bem}


\end{document}